\documentclass[12pt]{amsart}
\usepackage{fullpage,url}
\usepackage[dvips]{hyperref}
 \usepackage{amscd}   
\usepackage{amsmath}

\DeclareFontEncoding{OT2}{}{} 

\DeclareMathOperator{\lcm}{lcm}

\DeclareMathOperator{\ac}{ac}
\DeclareMathOperator{\tor}{tor}
\DeclareMathOperator{\Span}{Span}
\DeclareMathOperator{\alg}{alg}
\DeclareMathOperator{\hhat}{\hat{h}}

\DeclareMathOperator{\Frac}{Frac}

\newtheorem{theorem}{Theorem}[section]
\newtheorem{lemma}[theorem]{Lemma}

\theoremstyle{definition}
\newtheorem{definition}[theorem]{Definition}

\newtheorem{conjecture}[theorem]{Conjecture}
\newtheorem{statement}[theorem]{Statement}
\newtheorem{Claim}[theorem]{Claim}
\newtheorem{example}[theorem]{Example}

\theoremstyle{remark}
\newtheorem{remark}[theorem]{Remark}

\title{The local Lehmer inequality for Drinfeld modules}
\author{Dragos Ghioca}
\begin{document}
\begin{abstract}
We give a lower bound for the local height of a non-torsion element of a 
Drinfeld module.
\end{abstract}
\maketitle
\section{Introduction}
\footnotetext[1]{1991 AMS Subject Classification: Primary, 11G09; Secondary, 
11G50}
In this paper we will use the following notation: p is a prime number and $q$ is 
a power of $p$. We denote by $\mathbb{F}_q$ the finite field with $q$ elements. 
We let $C$ be a nonsingular projective curve defined over $\mathbb{F}_q$ and we 
fix a closed point on $C$, which we call $\infty$. Then we define $A$ as the 
ring of functions on $C$ that are regular everywhere except possibly at 
$\infty$. Also, in this paper, the elements of $\mathbb{F}_q^{\alg}$ will be 
called constants.

We let $K$ be a finitely generated field extension of $\mathbb{F}_q$. We fix a 
morphism $i:A\rightarrow K$. We define the operator $\tau$ as the power of the 
usual Frobenius with the property that for every $x\in K^{\alg}$, $\tau(x)=x^q$. 
Then we let $K\{\tau\}$ be the ring of polynomials in $\tau$ with coefficients 
from $K$.

A Drinfeld module is a morphism $\phi:A\rightarrow K\{\tau\}$ for which the 
coefficient of $\tau^0$ in $\phi_a$ is $i(a)$ for every $a\in A$. Following the 
definition from \cite{Goss} we will call $\phi$ a Drinfeld module of generic 
characteristic if $\ker(i)=\{0\}$ and we will call $\phi$ a Drinfeld module of 
finite characteristic if $\ker(i)\ne \{0\}$. In case of a Drinfeld module of 
generic characteristic we will identify $i(a)$ with $a$ for every $a\in A$.

In Section $2$ we will develop the theory of heights on Drinfeld modules. We 
denote by $\hhat:K^{\alg}\rightarrow \mathbb{R}_{\ge 0}$ the global height 
associated to a Drinfeld module $\phi$ and for each divisor $v$ (as defined in 
Section $2$) we define by $\hhat_v:K^{\alg}\rightarrow\mathbb{R}_{\ge 0}$ the 
corresponding local height associated to $\phi$.

The paper \cite{Den} proposed the following conjecture, which is the Drinfeld 
module analogue of the classical Lehmer inequality.
\begin{conjecture}
\label{C:Con1}
For the Drinfeld module $\phi:A\rightarrow K\{\tau\}$ there exists a constant 
$C>0$, depending only on $\phi$, such that any non-torsion point $x\in K^{\alg}$ 
satisfies $\hhat(x)\ge\frac{C}{[K(x):K]}$.
\end{conjecture}
A partial result towards this conjecture was obtained in \cite{DG}. 

The following statement would imply \eqref{C:Con1} and we refer to it as the 
local case of the Lehmer inequality for Drinfeld modules.
\begin{statement}
\label{C:Con2}
For the Drinfeld module $\phi:A\rightarrow K\{\tau\}$ there exists a constant 
$C>0$, depending only on $\phi$, such that for any $x\in K^{\alg}$ and any place 
$v$ of $K(x)$, if $\hhat_v(x)>0$, then $\hhat_v(x)\ge\frac{C}{[K(x):K]}$.
\end{statement}
In the third section of this paper we will prove that 
\eqref{C:Con2} is false but in the case of Drinfeld modules of finite 
characteristic there is the following result. 
\begin{theorem}
\label{T:T1'}
For $\phi :A\rightarrow K\{\tau\}$ a finite characteristic Drinfeld module, 
there exist two positive constants $C$ and $k$ depending only on $\phi$ such 
that if $x\in K^{\alg}$ and $v$ is a place of $K(x)$ for 
which $\hhat_v(x)>0$, then $\hhat_v(x)\ge\frac{C}{d^k}$ (where $d=[K(x):K]$).
\end{theorem}
Theorem \eqref{T:T1'} will follow from the following stronger result.
\begin{theorem}
\label{T:T1}
Let $\phi:A\rightarrow K\{\tau\}$ be a Drinfeld module of finite characteristic. 
Let $x\in K^{\alg}$ 
and let $d=[K(x):K]$. Let $v\in M_{K(x)}$ such that $\hhat_v(x)>0$. Denote by 
$v_0$ the place of $K$ sitting below $v$ and let $e(v|v_0)$ be the corresponding 
ramification index. 

There exists $C>0$ and $k\ge 1$, both depending only on $\phi$, such that 
$\hhat_v(x)\ge\frac{C}{e(v|v_0)^{k-1}d}$.
\end{theorem}
Moreover if $p$ does not divide $e(v\vert v_0)$, then we can give a very easy 
expression 
for the exponent $k$ in \eqref{T:T1} which will be optimal in 
this case as shown by example \eqref{E:E1}. This will be proved in theorem 
\eqref{T:T2}.

If $\phi$ is a Drinfeld module of generic characteristic, example \eqref{E:E2} 
will show that $\hhat_v(x)$ can be arbitrarily small and strictly positive 
regardless of $d=[K(x):K]$. In theorem \eqref{T:T3}, we will give the best 
result towards conjecture \eqref{C:Con2} for Drinfeld modules of generic 
characteristic.

We thank Bjorn Poonen and Thomas Scanlon for expositional suggestions. We 
express our gratitude to Thomas Scanlon for his encouragement and for asking 
the mathematical questions that led us to conjecture statement \eqref{C:Con2} 
which constituted the starting point for this paper.

\section{Heights associated to Drinfeld modules}

As stated in Section $1$, we are working with a Drinfeld module 
$\phi:A\rightarrow K\{\tau\}$. Because $K$ is a finitely generated field over 
$\mathbb{F}_q$, it is the 
function field of a variety $V$ defined over $\mathbb{F}_{q^m}$ for some 
$m\ge 1$ and in addition we can take $V$ to be normal and projective, embedded 
in $\mathbb{P}^M$ (for some $M\ge 1$). We define $M_K$ as the set of valuations 
of 
$K$ that are associated to irreducible divisors of $V$, i.e. codimension $1$ 
subvarieties of $V$. Then to an element 
$x$ from $K$, we associate its divisor 
$$(x)=\sum_{\rho\in M_K} v_{\rho}(x)\cdot\rho$$
where by $v_{\rho}(x)$ we denoted the order of $x$ at $\rho$.

For each $\rho\in M_K$, we denote by $\deg(\rho)$ the projective degree of 
$\rho$ 
in $\mathbb{P}^M$, which is the intersection number of $\rho$ with a generic 
$(M-N+1)$-dimensional hyperplane in $\mathbb{P}^M$ ($N=\dim\space V$). The 
following product formula holds  $$\sum_{\rho\in 
M_K}\deg(\rho)\cdot v_{\rho}(x)=0\text{.}$$ 
For simplicity of notation in the rest of this paper we will drop the index 
$\rho$ from the valuation $v$.

Now we construct the local heights $\hhat_v$ with respect to the Drinfeld 
module $\phi$. Our construction follows \cite{Poo} together with the 
observations 
from \cite{Wan} that extend the construction to finitely generated function 
fields. So, for $x\in K$ and $v\in M_K$, we set $\tilde{v}(x)=\min\{ 0,v(x)\}$ 
and for a nonconstant element $a\in A$, we define 
$$V_v(x)=\lim_{n\rightarrow\infty}\frac{\tilde{v}(\phi_{a^n}(x))}{\deg 
(\phi_{a^n})}.$$
This function satisfies the same properties as in Propositions $1$-$4$ from 
\cite{Poo}. We define $$\hhat_v(x)=-\deg (v) V_v(x)$$
where $\deg (v)$ is the degree of the divisor $v$ as defined above.

This defines the local 
heights only for elements of $K$, but we will be interested in 
extending 
them to $K^{\alg}$. For this, let $x\in L$, where $L$ is a finite extension of 
$K$. We let $W$ be the normalization of $V$ in $L$ and form the set $M_L$ of 
valuations of $L$ associated to $W$. As shown in \cite{Wan}, for 
every $v\in M_K$, there exist finitely many $w\in M_L$ extending $v$. When we 
work in 
such a setting, our convention will always be that the valuations are functions 
with range $\mathbb{Z}$. Thus $w\vert_K=e(w\vert v)v$, where $e(w\vert v)$ is 
the corresponding ramification index. We define the function $V_w$ as $V_v$ from 
above. Then we let
$$\hhat_w(x)=-\frac{\deg(v)f(w\vert v)}{[L:K]}V_w(x)$$ where $f(w\vert v)$ is 
the 
relative residue degree between the residue field of $L$ at $w$ and the residue 
field of $K$ at $v$. 

Then, just as in \cite{Poo}, we define the global height with respect to 
$\phi$ by
$$\hhat(x)=\sum_{w\in M_L}\hhat_w(x).$$
The above sum is finite due to a similar argument as the one from Proposition 
$6$ of \cite{Poo}. 

If $L'$ is a finite extension of $L$ and $w'|w$ is any valuation on $L'$ 
extending $w$, then $V_{w'}(x)=e(w'|w)V_w(x)$. Because 
$\sum_{w'|w}e(w'|w)f(w'|w)=[L':L]$, we get $\sum_{w'|w}\hhat_{w'}(x)=\hhat_w(x)$ 
for every $x\in L$ and every $w\in M_L$. Thus, our definition of the global 
height is independent of the field $L$ containing $x$.

Let $t$ be a non-constant element of $A$ and 
$\phi_t=\sum_{i=r_0}^{r}a_i\tau^{i}$, with $a_{r_0}\ne 0$. The first results of 
this section do not rely on $\phi$ being of finite characteristic or not and so, 
we do not specify right now if $r_0=0$ or $r_0>0$.

In proving \eqref{T:T1} we may replace $K$ by a finite extension $K'$. This will 
only induce a constant factor  $[K':K]$ in the denominator of the 
lower bound for the local height. Also, \eqref{T:T1} is not affected if we 
replace $\phi$ by a Drinfeld module that is isomorphic to $\phi$. 
Thus we can conjugate $\phi$ by an element $\gamma\in K^{\alg}$ such that 
$\phi^{(\gamma )}$, the conjugated Drinfeld module, has the property that 
$\phi^{(\gamma )}_t$ is monic as a polynomial in $\tau$. Then $\phi$ and 
$\phi^{(\gamma )}$ are isomorphic over $K(\gamma)$ which is a finite extension 
of $K$ (because $\gamma$ satisfies the equation $\gamma^{q^r-1}a_r=1$).

So, we will prove theorem \eqref{T:T1} for $\phi^{(\gamma)}$ and because 
$\hhat_{\phi ,v}(x)=\hhat_{\phi^{(\gamma)},v}(\gamma^{-1}x)$ (as proved in 
\cite{Poo}, Proposition $2$) the result will follow for $\phi$. For simplicity 
of notation 
we will suppose from now on that $\phi_t$ is monic as a polynomial in $\tau$.

Let $x$ be a nonzero element of $K^{\alg}$ and let $L=K(x)$. Denote by $S$ the 
finite subset of $M_L$ where the coefficients $a_i$, for 
$i\in\{r_0,\dots,r-1\}$, have 
poles. Also, 
denote by $S_0$ the finite set of divisors from $M_K$ where the coefficients 
$a_i$ 
have poles. Thus, each divisor from $S$ 
sits above an unique divisor from $S_0$.

For each $v\in M_L$ denote by 
\begin{equation}
\label{E:defM_v}
M_v=\min_{i\in\{r_0,\dots,r-1\}}\frac{v(a_i)}{q^r-q^i}
\end{equation}
where by convention: $v(0)=+\infty$. If $r_0=r$, definition \eqref{E:defM_v} is 
void and in that case we define $M_v=+\infty$.

Note that $M_v<0$ if and only if $v\in S$.

For each $v\in S$ we fix a uniformizer $\pi_v\in L$ of the place $v$. We define 
next the concept of angular component for every $y\in L\setminus\{0\}$.

\begin{definition}
\label{D:angular}
Assume $v\in S$. For every nonzero $y\in L$ we define the angular component of 
$y$ at $v$, denoted by $\ac_{\pi_v}(y)$, the residue at $v$ of $y\pi_v^{-v(y)}$. 
(Note that the angular component is never $0$.)
\end{definition}
We can define in a similar manner as above the notion of angular component at 
each $v\in M_L$ but we will work with angular components at the places from $S$ 
only.

The main property of the angular component is that for every $y,z\in L$, 
$v(y-z)>v(y)=v(z)$ if and only if $(v(y),\ac_{\pi_v}(y))=(v(z),\ac_{\pi_v}(z))$. 

Our strategy for proving \eqref{T:T1} will be to prove that if $\hhat_v(x)>0$ 
then \emph{either} 
$$\hhat_v(x)\ge\frac{C}{e(v|v_0)^{\frac{r}{r_0}-1}d}$$
where $e(v|v_0)$ is the corresponding ramification index of $v$ over $K$, 
$d=[L:K]$ and $C>0$ is a constant depending only on $\phi$, \emph{or}
$$v\in S\text{ and }(v(x),\ac_{\pi_v}(x))\text{ belongs to a set of cardinality 
we can control.}$$

For $v\in S$ we define 
\begin{equation}
\label{E:defP_v}
P_v=\left\{\frac{v(a_i)-v(a_j)}{q^j-q^i}|r_0\le i<j\le r\text{ and } a_i\ne 0\ne 
a_j\right\}\cup\{0\}.
\end{equation}
Clearly, $|P_v|\le 1+\binom{r-r_0+1}{2}$. For each $\alpha\in P_v$ we let $l\ge 
1$ and let $ i_0<i_1<\dots <i_l$ be all the indices $i$ for which $a_i\ne 0$ and 
moreover, for $j,k\in\{0,\dots,l\}$ with $j\ne k$, we have
\begin{equation}
\label{E:Newton}
\frac{v(a_{i_j})-v(a_{i_k})}{q^{i_k}-q^{i_j}}=\alpha .
\end{equation}
We define $R_v(\alpha)$ as the set containing $\{1\}$ and all the nonzero 
solutions of the equation 
\begin{equation}
\label{E:defR_v}
\sum_{j=0}^{l}\ac_{\pi_v}(a_{i_j})X^{q^{i_j}}=0.
\end{equation}
Clearly, for every $\alpha\in P_v$, $|R_v(\alpha)|\le q^r$.

Note that if $\alpha=0$, there might be no indices $i_j$ and $i_k$ as in 
\eqref{E:Newton}. In that case, the construction of $R_v(0)$ from 
\eqref{E:defR_v} is void and so, we define $R_v(0)=\{1\}$. The motivation for 
the special case $0\in P_v$ and $1\in R_v(0)$ is explained in the proof of lemma 
\eqref{L:L5'}.

\begin{lemma}
\label{L:L0}
Assume $v\in S$. If $v(\phi_t(x))>\min_{i\in\{r_0,\dots,r\}}v(a_ix^{q^i})$ then 
$(v(x),\ac_{\pi_v}(x))\in P_v\times R_v(v(x))$.
\end{lemma}
\begin{proof}
If $v(x)>\min_{i\in\{r_0,\dots,r\}}v(a_ix^{q^i})$ it means that there exists 
$l\ge 1$ and $ i_0<\dots <i_l$ such that
\begin{equation}
\label{E:0,l,v}
v(a_{i_0}x^{q^{i_0}})=\dots =v(a_{i_l}x^{q^{i_l}})
\end{equation} 
and also 
\begin{equation}
\label{E:0,l,ac}
\sum_{j=0}^{l}\ac_{\pi_v}(a_{i_j})\ac_{\pi_v}(x)^{q^{i_j}}=0.
\end{equation}
Equations \eqref{E:0,l,v} and \eqref{E:0,l,ac} yield $v(x)\in P_v$ and 
$\ac_{\pi_v}(x)\in R_v(v(x))$ respectively, according to \eqref{E:defP_v} and 
\eqref{E:defR_v}.
\end{proof}

\begin{lemma}
\label{L:L2'}
Let $v\in M_L$ and let $v_0\in M_K$ be the unique 
valuation of $K$ sitting below $v$. If $v(x)<\min\{0,M_v\}$, then 
$\hhat_v(x)=\frac{-\deg(v_0)f(v\vert 
v_0)}{[L:K]}v(x)$.
\end{lemma}
\begin{proof}  
For every $i\in\{r_0,\dots,r-1\}$, $v(a_ix^{q^i})=v(a_i)+q^iv(x)>q^rv(x)$ 
because $v(x)<M_v=\min_{i\in\{r_0,\dots,r-1\}}\frac{v(a_i)}{q^r-q^i}$. This 
shows that $v(\phi_t(x))=q^rv(x)<v(x)<\min\{0,M_v\}$. By induction, 
$v(\phi_{t^n}(x))=q^{rn}v(x)$ 
for all $n\ge 1$. So, $V_v(x)=v(x)$ 
and 
$$\hhat_v(x)=\frac{-\deg(v_0)f(v\vert v_0)}{[L:K]}v(x).$$
\end{proof}

An immediate corollary to \eqref{L:L2'} is the following result.
\begin{lemma}
\label{L:L1'}
Assume $v\notin S$. If $v(x)<0$ then $\hhat_v(x)=\frac{-\deg(v_0)f(v\vert 
v_0)}{[L:K]}v(x)$, 
while if $v(x)\ge 0$ then $\hhat_v(x)=0$.
\end{lemma}
\begin{proof}  
First, it is clear that if $v(x)\ge 0$ then for all $ n\ge 1$, 
$v(\phi_{t^n}(x))\ge 0$ because all the coefficients of $\phi_t$ and thus of 
$\phi_{t^n}$ have non-negative valuation at $v$. Thus, $V_v(x)=0$ and so 
$$\hhat_v(x)=0.$$
Now, if $v(x)<0$, then $v(x)<M_v$ because $M_v\ge 0$ ($v\notin S$). So, applying 
the result of \eqref{L:L2'} we conclude the proof of this lemma.
\end{proof}

We will get a better insight into the local heights behaviour with the following 
lemma.
\begin{lemma}
\label{L:L3'}
Assume $v\in S$ and $v(x)\le 0$. If $(v(x),\ac_{\pi_v}(x))\notin P_v\times 
R_v(v(x))$ then $v(\phi_t(x))<M_v$, unless $q=2$, $r=1$ and $v(x)=0$.
\end{lemma}
\begin{proof}
Lemma \eqref{L:L0} implies that there exists $i_0\in\{r_0,\dots,r\}$ such that 
for all $i\in\{r_0,\dots,r\}$ we have 
$v(a_ix^{q^i})\ge v(a_{i_0}x^{q^{i_0}})=v(\phi_t(x))$. 

Suppose \eqref{L:L3'} is not true and so, there exists $j_0<r$ such that
$$\frac{v(a_{j_0})}{q^r-q^{j_0}}\le v(\phi_t(x))=v(a_{i_0})+q^{i_0}v(x).$$
This means that
\begin{equation}
\label{E:a_j_0_1}
v(a_{j_0})\le (q^r-q^{j_0})v(a_{i_0})+(q^{r+i_0}-q^{i_0+j_0})v(x).
\end{equation}
On the other hand, by our assumption about $i_0$, we know that 
$v(a_{j_0}x^{q^{j_0}})\ge v(a_{i_0}x^{q^{i_0}})$ which means that 
\begin{equation}
\label{E:a_j_0_2}
v(a_{j_0})\ge v(a_{i_0})+(q^{i_0}-q^{j_0})v(x).
\end{equation}
Putting together inequalities \eqref{E:a_j_0_1} and \eqref{E:a_j_0_2}, we get 
$$v(a_{i_0})+(q^{i_0}-q^{j_0})v(x)\le 
(q^r-q^{j_0})v(a_{i_0})+(q^{r+i_0}-q^{i_0+j_0})v(x).$$
Thus 
\begin{equation}
\label{E:1'}
v(x)(q^{r+i_0}-q^{i_0+j_0}-q^{i_0}+q^{j_0})\ge -v(a_{i_0})(q^r-q^{j_0}-1).
\end{equation}
But
$q^{r+i_0}-q^{i_0+j_0}-q^{i_0}+q^{j_0}=q^{r+i_0}(1-q^{j_0-r}-q^{-r}+q^{j_0-r-i_0
})$ and because $j_0<r$ and $q^{j_0-r-i_0}>0$, we obtain
\begin{equation}
\label{E:2'}
1-q^{j_0-r}-q^{-r}+q^{j_0-r-i_0}>1-q^{-1}-q^{-r}\ge 1-2q^{-1}\ge 0.
\end{equation}
Also, $q^r-q^{j_0}-1\ge q^r-q^{r-1}-1=q^{r-1}(q-1)-1\ge 0$ with equality if and 
only if 
$q=2$, $r=1$ and $j_0=0$. We will analyze this case separately. So, as long as 
we are not in this special case, we do have 
\begin{equation}
\label{E:3'}
q^r-q^{j_0}-1>0.
\end{equation}
Now we have two possibilities (remember that $v(x)\le 0$):

(i)  $v(x)<0$  

In this case, \eqref{E:1'}, \eqref{E:2'} and \eqref{E:3'} tell us that 
$-v(a_{i_0})<0$. Thus, $v(a_{i_0})>0$. But we know from our hypothesis on 
$i_0$ that $v(a_{i_0}x^{q^{i_0}})\le v(x^{q^r})$ which is in contradiction with 
the combination of the following facts: $v(x)<0$, $i_0\le r$ and $v(a_{i_0})>0$.

(ii)  $v(x)=0$

Then another use of \eqref{E:1'}, \eqref{E:2'} and \eqref{E:3'} gives us 
$-v(a_{i_0})\le 0$; thus $v(a_{i_0})\ge 0$. This would mean that 
$v(a_{i_0}x^{q^{i_0}})\ge 0$ and this contradicts our choice for $i_0$ because 
we know from the fact that $v\in S$, that there exists $i\in \{r_0,\dots,r\}$ 
such that $v(a_i)<0$. So, then we would have 
$$v(a_ix^{q^i})=v(a_i)<0\le v(a_{i_0}x^{q^{i_0}}).$$
Thus, in either case (i) or (ii) we get a contradiction that proves the lemma 
except in the 
special case that we excluded above: $q=2$, $r=1$ and $j_0=0$. If we have $q=2$ 
and $r=1$ then 
$$\phi_t(x)=a_0x+x^2.$$
Note that if $a_0\in \mathbb{F}_p^{\alg}$, $S$ is empty and so, the result of 
our lemma is vacuously true. Thus, we suppose from now on that in this case: 
$q=2$ and $r=1$, $a_0\notin \mathbb{F}_p^{\alg}$ and so, $S$ consists of the 
irreducible divisors of the pole of $a_0$.

If $v(x)\le 0$, then either $v(x)<M_v=v(a_0)$, in which case again 
$v(\phi_t(x))<M_v$ (as shown in the proof of lemma \eqref{L:L2'}), or $v(x)\ge 
M_v$ and so, $i_0=0$ (because in this 
case $v(a_0x)\le v(x^2)$). In the latter case, 
$$v(\phi_t(x))=v(a_0x)=v(a_0)+v(x)<v(a_0)=M_v$$
unless $v(x)=0$. So, we see that indeed, only $v(x)=0$, $q=2$ and $r=1$ can make 
$v(\phi_t(x))\ge M_v$ in the hypothesis of \eqref{L:L3'}. 
\end{proof}

\begin{lemma}
\label{L:L4'}
Assume $v\in S$. Excluding the case $q=2$, $r=1$ and $v(x)=0$, we have that if 
$v(x)\le 0$ 
then either $\hhat_v(x)>\frac{-\deg(v_0)f(v\vert v_0)M_v}{q^r[L:K]}$ or 
$(v(x),ac_{\pi_v}(x))\in P_v\times R_v(v(x))$.
\end{lemma}

\begin{proof}
If $v(x)\le 0$ then

$$\text{\emph{either}: (i) } v(\phi_t(x))<M_v\text{ ,}$$

in which case by \eqref{L:L2'} we have that 
$\hhat_v(\phi_t(x))=\frac{-\deg(v_0)f(v\vert 
v_0)}{[L:K]}v(\phi_t(x))$. So, case (i) yields
\begin{equation}
\label{E:4'}
\hhat_v(x)=\frac{-\deg(v_0)f(v\vert 
v_0)}{[L:K]}\cdot\frac{v(\phi_t(x))}{\deg\phi_t} > 
\frac{-\deg(v_0)f(v\vert v_0)}{[L:K]}\cdot\frac{M_v}{q^r}
\end{equation}

$$\text{\emph{or}: (ii) }v(\phi_t(x))\ge M_v\text{ ,}$$

in which case, lemma \eqref{L:L3'} yields 
\begin{equation}
\label{E:(ii)}
v(\phi_t(x))>v(a_{i_0}x^{q^{i_0}})=\min_{i\in\{r_0\dots,r\}}v(a_ix^{q^i}).
\end{equation}
Using \eqref{E:(ii)} and lemma \eqref{L:L0} we conclude that case (ii) yields 
$(v(x),\ac_{\pi_v}(x))\in P_v\times R_v(v(x))$.
\end{proof}

Now we analyze the excluded case from lemma \eqref{L:L4'}.
\begin{lemma}
\label{L:L5'}
Assume $v\in S$. If $v(x)\le 0$ then either $(v(x),ac_{\pi_v}(x))\in P_v\times 
R_v(v(x))$ or $\hhat_v(x)\ge\frac{-\deg(v_0)f(v\vert v_0)M_v}{q^r[L:K]}$.
\end{lemma}

\begin{proof}
Using the result of \eqref{L:L4'} we have left to analyze the case: $q=2$, $r=1$ 
and $v(x)=0$. 

As shown in the proof of \eqref{L:L3'}, in this case $\phi_t(x)=a_0x+x^2$ and 
$$v(\phi_t(x))=v(a_0)=M_v<0.$$ 
Then, either $v(\phi_{t^2}(x))=v(\phi_t(x)^2)=2M_v<M_v$ or 
$v(\phi_{t^2}(x))>v(a_0\phi_t(x))=v(\phi_t(x)^2)$. If the former case holds, 
then by \eqref{L:L2'},
$$\hhat_v(\phi_{t^2}(x))=\frac{-\deg(v_0)f(v\vert v_0)}{[L:K]}\cdot 
2M_v\Rightarrow 
\hhat_v(x)=\frac{-\deg(v_0)f(v\vert v_0)}{[L:K]}\frac{2M_v}{4}.$$
If the latter case holds, i.e. 
$v(\phi_{t}(\phi_t(x)))>v(a_0\phi_t(x))=v(\phi_t(x)^2)$, it means that 
$\ac_{\pi_v}(\phi_t(x))$ satisfies the equation
$$\ac_{\pi_v}(a_0)X+X^2=0.$$ 
Because the angular component is never $0$, it must be that 
$\ac_{\pi_v}(\phi_t(x))=\ac_{\pi_v}(a_0)$ (remember that we are working now in 
characteristic $2$). But, because 
$v(a_0x)<v(x^2)$ we can relate the angular component of $x$ and the angular 
component of $\phi_t(x)$ and so,
$$ac_{\pi_v}(a_0)=\ac_{\pi_v}(\phi_t(x))=\ac_{\pi_v}(a_0x)= 
ac_{\pi_v}(a_0)\ac_{\pi_v}(x).$$
This means $\ac_{\pi_v}(x)=1$ and so, the excluded case amounts to a 
dichotomy similar to the one from \eqref{L:L4'}: either 
$(v(x),\ac_{\pi_v}(x))=(0,1)$ or $\hhat_v(x)=\frac{\deg(v_0)f(v\vert 
v_0)}{[L:K]}\frac{-M_v}{2}$. The definitions of $P_v$ and $R_v(\alpha)$ from 
\eqref{E:defP_v} and \eqref{E:defR_v} respectively, yield that $(0,1)\in 
P_v\times R_v(0)$.
\end{proof}

Finally, we note that in \eqref{L:L5'} we have 
$$-\deg(v_0)f(v\vert v_0)M_v\ge -M_v> \frac{e(v\vert v_0)}{q^r}.$$
We have obtained the following dichotomy.
\begin{lemma}
\label{L:L6'}
Assume $v\in S$. If $v(x)\le 0$ then either 
$\hhat_v(x)\ge\frac{e(v|v_0)}{q^{2r}d}$ or $(v(x),\ac_{\pi_v}(x))\in P_v\times 
R_v(v(x))$ with $|P_v|\le 1+\binom{r-r_0+1}{2}$ and for each $\alpha\in P_v$, 
$|R_v(\alpha)|\le q^r$.
\end{lemma}

\begin{lemma}
\label{L: different values}
There are no $x$ and $x'$ verifying the following properties 

(a) $v(x)\ne v(x')$; 

(b) $(v(x),\ac_{\pi_v}(x))\notin P_v\times R_v(v(x))$ and 
$(v(x'),\ac_{\pi_v}(x'))\notin 
P_v\times R_v(v(x'))$;

(c) $v(\phi_t(x))=v(\phi_t(x'))$.
\end{lemma}
\begin{proof}  
Suppose \eqref{L: different values} is not true and so, there exist  $x$, $x'$ 
satisfying (a),(b) and (c). Property (b) and lemma \eqref{L:L0} yield that there 
exists $ i_1\in\{r_0,\dots,r\}$ such that
\begin{equation}
\label{E:4}
v(\phi_t(x))=v(a_{i_1})+q^{i_1}v(x)
\end{equation}
and for every $i\in\{r_0,\dots,r\}$,
\begin{equation}
\label{E:5}
v(a_{i_1})+q^{i_1}v(x)\le v(a_i)+q^iv(x).
\end{equation}
We have the similar equations for $x'$ and for some $i_2\in\{r_0,\dots,r\}$,
\begin{equation}
\label{E:6}
v(\phi_t(x'))=v(a_{i_2})+q^{i_2}v(x')
\end{equation}
where for every $i\in\{r_0,\dots,r\}$,
\begin{equation}
\label{E:7}
v(a_{i_2})+q^{i_2}v(x')\le v(a_i)+q^iv(x').
\end{equation}
Also, we supposed that we have 
\begin{equation}
\label{E:8}
v(\phi_t(x))=v(\phi_t(x')).
\end{equation}
We assume that $i_1\ne i_2$, because $i_1=i_2$ would imply from \eqref{E:4}, 
\eqref{E:6} and \eqref{E:8} that $v(x)=v(x')$. So, without loss of generality we 
may assume that $i_1<i_2$. We use \eqref{E:5} for 
$i=i_2$ and so, we get $v(a_{i_1})+q^{i_1}v(x)\le v(a_{i_2})+q^{i_2}v(x)$ which 
implies
\begin{equation}
\label{E:9}
v(x)\ge\frac{v(a_{i_1})-v(a_{i_2})}{q^{i_2}-q^{i_1}}.
\end{equation}
Now, using \eqref{E:7} with $i=i_1$, we get $v(a_{i_2})+q^{i_2}v(x')\le 
v(a_{i_1})+q^{i_1}v(x')$ which implies 
\begin{equation}
\label{E:10}
v(x')\le\frac{v(a_{i_1})-v(a_{i_2})}{q^{i_2}-q^{i_1}}.
\end{equation}
But because of \eqref{E:8} together with \eqref{E:4} and \eqref{E:6}, we have 
$v(a_{i_1})+q^{i_1}v(x)=v(a_{i_2})+q^{i_2}v(x')$ which implies that 
\begin{equation}
\label{E:11}
v(x')=\frac{q^{i_1}v(x)+v(a_{i_1})-v(a_{i_2})}{q^{i_2}}.
\end{equation}
Using \eqref{E:10} and \eqref{E:11}, we get 
$$\frac{q^{i_1}v(x)+v(a_{i_1})-v(a_{i_2})}{q^{i_2}}\le\frac{v(a_{i_1})-v(a_ 
{i_2})}{q^{i_2}-q^{i_1}}.$$
So, $q^{i_1}(q^{i_2}-q^{i_1})v(x)\le q^{i_1}(v(a_{i_1})-v(a_{i_2}))$, which 
implies that 
$$v(x)\le \frac{v(a_{i_1})-v(a_{i_2})}{q^{i_2}-q^{i_1}}$$ 
which combined with \eqref{E:9} shows that 
$$v(x)=\frac{v(a_{i_1})-v(a_{i_2})}{q^{i_2}-q^{i_1}}.$$
Then, using \eqref{E:11} we get 
$$v(x')=\frac{q^{i_1} \frac{v(a_{i_1})-v(a_{i_2})}{q^{i_2}-q^{i_1}} 
+v(a_{i_1})-v(a_{i_2})}{q^{i_2}}=\frac{v(a_{i_1})-v(a_{i_2})}{q^{i_2}-q^{i_1}}=v
(x)$$
which shows that indeed $v(\phi_t(x))$ can be obtained from an unique value for 
$v(x)$.
\end{proof}

\begin{lemma}
\label{L:Lacv}
Assume $v\in S$. If $(v(x),\ac_{\pi_v}(x))\notin P_v\times R_v(v(x))$ then for 
each of the values $(\alpha_1,\gamma_1)=(v(\phi_t(x)),\ac_{\pi_v}(\phi_t(x)))$ 
there are at most $q^r$ possible values $\gamma$ for $\ac_{\pi_v}(x)$.
\end{lemma}

\begin{proof}
Indeed, we saw in lemma \eqref{L: different values} that $v(x)$ is uniquely 
determined given $\alpha_1=v(\phi_t(x))$ under the hypothesis of \eqref{L:Lacv}. 
We also have
\begin{equation}
\label{E:12}
\ac_{\pi_v}(\phi_t(x))=\sum_{j} \ac_{\pi_v}(a_{i_j})\ac_{\pi_v}(x)^{q^{i_j}}
\end{equation}
where $i_j$ runs through a prescribed subset of $\{r_0,\dots,r\}$ corresponding 
to those $i$ such that $v(a_i)+q^iv(x)=v(\phi_t(x))$. 
This subset of indices $i_j$, depends only on $\alpha_1=v(x)$. So, there are at 
most $q^r$ possible values for $\ac_{\pi_v}(x)$ to solve \eqref{E:12} given 
$\gamma_1=\ac_{\pi_v}(\phi_t(x))$.
\end{proof}

From now on in this section we will suppose that

\fbox{$r_0\ge 1$, i.e. $\phi$ has finite 
characteristic and $t\in A$ has the property that $\phi_t$ is inseparable.}

\begin{lemma}
\label{L:L11}
For $v\in S$ denote by $N_v=\max \left\{\frac{-v\left(a_i\right)}{q^i-1} \mid 
1\le i\le r \right\} $ (remember our convention on $v(0)=+\infty$). If 
$v\left(x\right)\ge N_v$, then $\hhat_v(x)=0$.
\end{lemma} 
\begin{proof}

Indeed, if $v(x)\ge N_v$ then $v\left(\phi_t\left(x\right)\right)\ge \min_{1\le 
i \le r}\lbrace 
q^iv\left(x\right)+v\left(a_i\right)\rbrace\ge v\left(x\right)\ge N_v$. By 
induction, we get that $v(\phi_{t^n}(x))\ge N_v$ for all $n\ge 1$, which yields 
that $V_v(x)=0$ and so,
$$\hat{h}_v\left(x\right)=0.$$
\end{proof} 

Thus, if $v\in S$ and $\hhat_v(x)>0$ it must be that $v\left(x\right)<N_v$.

\begin{lemma}
\label{L:L1}
Assume $v\in S$. If $v(x)<N_v$ and if $(v(x),\ac_{\pi_v}(x))\notin P_v\times 
R_v(v(x))$ then 
$v(\phi_t(x))<v(x)$.
\end{lemma}
\begin{proof}  
Indeed, by the hypothesis and by lemma \eqref{L:L0}, there exists 
$i_0\in\{r_0,\dots,r\}$ such that 
for all $ i\in \{r_0,\dots,r\}$,
\begin{equation}
\label{E:13}
v(a_{i_0})+q^{i_0}v(x)=v(\phi_t(x))\le v(a_i)+q^iv(x).
\end{equation}
If $v(\phi_t(x))\ge v(x)$ then, using \eqref{E:13}, we get that
$$v(x)\le v(a_i)+q^iv(x)$$ 
which implies that $v(x)\ge -\frac{v(a_i)}{q^i-1}$ for every $i$. Thus
$$v(x)\ge N_v,$$
contradicting the hypothesis of our lemma. So, we must have $v(\phi_t(x))<v(x)$. 
In particular, we also get that 
$v(a_{i_0})+q^{i_0}v(x)<v(x)$, i.e.
\begin{equation}
\label{E:vi_0}
v(x)<\frac{-v(a_{i_0})}{q^{i_0}-1}.
\end{equation}
\end{proof}

Our goal is establishing a dichotomy similar to the one from lemma \eqref{L:L6'} 
under the following hypothesis

\fbox{$v\in S$, $\hhat_v(x)>0$ and $0<v(x)<N_v$.}

In lemma \eqref{L:L1} we saw that if $v(x)<N_v$ then either 
$(v(x),\ac_{\pi_v}(x))\in P_v\times R_v(v(x))$ or $v(\phi_t(x))<v(x)$. In the 
latter case, if $v(\phi_t(x))>0$ we apply then the same reasoning to $\phi_t(x)$ 
and derive that either $(v(\phi_t(x)),\ac_{\pi_v}(\phi_t(x)))\in P_v\times 
R_v(v(\phi_t(x)))$ or $v(\phi_{t^2})<v(\phi_t(x))$. We repeat this analysis and 
after a finite number of steps, say $n$, we must have that either 
$v(\phi_{t^n}(x))\le 0$ or $(v(\phi_{t^n}(x)),\ac_{\pi_v}(\phi_{t^n}(x)))\in 
P_v\times R_v(v(\phi_{t^n}(x)))$.  But we analyzed in \eqref{L:L6'} what happens 
to the cases in which, for an element $y$ of positive local height at $v$, 
$v(y)\le 0$. We obtained that either
\begin{equation}
\label{E:14}
\hat{h}_v(y)\ge \frac{e(v\vert v_0)}{q^{2r}d}
\end{equation}
or
\begin{equation}
\label{E:15}
(v(y),\ac_{\pi_v}(y))\in P_v\times R_v(v(y))
\end{equation}
and $|P_v|\le 1+\binom{r-r_0+1}{2}\le 1+\frac{r^2-r}{2}=\frac{r^2-
r+2}{2}$ because $r_0\ge 1$.

We will use repeatedly equations \eqref{E:14} and \eqref{E:15} for 
$y=\phi_{t^n}(x)$. So, if 
\eqref{E:14} holds for $y=\phi_{t^n}(x)$ then 
\begin{equation}
\label{E:16}
\hat{h}_v(x)\ge\frac{e(v\vert v_0)}{q^{rn}q^{2r}d}.
\end{equation}
We will see next what happens if \eqref{E:15} holds. We can go back through the 
steps that we made in order to get to \eqref{E:15} and see that actually $v(x)$ 
and $\ac_{\pi_v}(x)$ belong to prescribed sets of cardinality independent of 
$n$.

\begin{lemma}
\label{L:L2}
Assume $v\in S$ and suppose that $v(x)<N_v$. If 
$(v(\phi_{t^k}(x)),\ac_{\pi_v}(\phi_{t^k}(x)))\notin P_v\times 
R_v(v(\phi_{t^k}(x)))$ for $0\le k\le n-1$, then for each value 
$(\alpha_n,\gamma_n)=(v(\phi_{t^n}(x)),\ac_{\pi_v}(\phi_{t^n}(x)))$, $v(x)$ is 
uniquely determined and $\ac_{\pi_v}(x)$ belongs to a set of cardinality at most 
$q^{r\cdot\binom{r}{2}}$.
\end{lemma}
\begin{proof} 
The fact that $v(x)$ is uniquely determined follows after $n$ successive 
applications of lemma \eqref{L: different values} to $\phi_{t^{n-1}}(x),\dots, 
\phi_t(x),x$.

Because $(v(\phi_{t^k}(x)),\ac_{\pi_v}(\phi_{t^k}(x)))\notin P_v\times 
R_v(v(\phi_{t^k}(x)))$ for $k<n$, it means that we are solving an equation of 
the form 
\begin{equation}
\label{E:17}
\sum_j \ac_{\pi_v}(a_{i_j}) 
\ac_{\pi_v}(\phi_{t^k}(x))^{q^{i_j}}=\ac_{\pi_v}(\phi_{t^{k+1}}(x))
\end{equation}
in order to express $\ac_{\pi_v}(\phi_{t^k}(x))$ in terms of 
$\ac_{\pi_v}(\phi_{t^{k+1}}(x))$ for each $k<n$. The equations \eqref{E:17} are 
uniquely determined by the sets of indices $i_j\in\{r_0,\dots,r\}$ which in 
turn are uniquely determined by $v(\phi_{t^k}(x))$, i.e. for each $k$ and each 
corresponding index $i_j$
\begin{equation}
\label{E:vdependence}
v(a_{i_j}\phi_{t^k}(x)^{q^{i_j}})=\min_{i\in\{r_0,\dots,r\}}v(a_i\phi_{t^k}(x)^{
q^i}).
\end{equation}

Using the result of 
\eqref{L:L1} and the hypothesis of our lemma, we see that
\begin{equation}
\label{E:18}
v(x)>v(\phi_t(x))>v(\phi_{t^2}(x))>\dots >v(\phi_{t^n}(x))
\end{equation}
and so the equations from \eqref{E:17} appear in a prescribed order. Now, in 
most of the cases, these equations will consist of only one term on their 
left-hand side; i.e. they will look like 
\begin{equation}
\label{E:19}
\ac_{\pi_v}(a_{i_0})\ac_{\pi_v}(\phi_{t^k}(x))^{q^{i_0}}=\ac_{\pi_v} 
(\phi_{t^{k+1}}(x)).
\end{equation}
Equation \eqref{E:19} has an unique solution. The other equations of type 
\eqref{E:17} but not of type \eqref{E:19} are associated to some of the values 
of $v(\phi_{t^k}(x))\in P_v$. Indeed, according to the definition of $P_v$ from 
\eqref{E:defP_v}, only for those values we can have for $i\ne i'$ 
\begin{equation}
\label{E:number}
v(a_i)+q^iv(x)=v(a_{i'})+q^{i'}v(x)
\end{equation}
and so, both indices $i$ and $i'$ can appear in \eqref{E:17}.

So, the number of equations of type \eqref{E:17} but not of type \eqref{E:19} is 
at most $\binom{r}{2}$ (remember that we are working under the assumption that 
$\phi_t$ is inseparable, i.e. $r_0\ge 1$). Moreover these equations will appear 
in a 
prescribed order, each not more than once, because of \eqref{E:18}. These 
observations determine the construction of the finite set that will contain all 
the possible values for $\ac_{\pi_v}(x)$, given 
$\gamma_n=\ac_{\pi_v}(\phi_{t^n}(x)))$. An equation of type \eqref{E:17} can 
have at most 
$q^r$ solutions; thus $\ac_{\pi_v}(x)$ lives in a set of cardinality at most 
$q^{r\cdot\binom{r}{2}}$.
\end{proof}

Because of the result of \eqref{L:L2}, we know that we can construct in an 
unique way $v(x)$ given $v(\phi_{t^n}(x))$ and the fact that for every $j<n$, 
$\phi_{t^j}(x)$ does not satisfy \eqref{E:15}. So, for each $n$ there are at 
most $\vert P_v\vert$ values for $v(x)$ such that 
\begin{equation}
\label{E:nmin}
(v(\phi_{t^n}(x)),\ac_{\pi_v}(\phi_{t^n}(x)))\in P_v\times R_v(v(\phi_{t^n}(x)))
\end{equation}
where $n$ is minimal with this property. We denote by $P_v(n)$ this set of 
values for $v(x)$. By convention: $P_v(0)=P_v$.

Also, lemma \eqref{L:L2} yields that for each fixed 
$(v(\phi_{t^n}(x)),\ac_{\pi_v}(\phi_{t^n}(x)))\in P_v\times 
R_v(v(\phi_{t^n}(x)))$ there are at most 
\begin{equation}
\label{E:possibilities}
q^{r\cdot\binom{r}{2}}=q^{\frac{r^3-r^2}{2}}
\end{equation}
possibilities for $\ac_{\pi_v}(x)$. For $\alpha=v(x)\in P_v(n)$ we define by 
$R_v(\alpha)$ the set of all possible values for $\ac_{\pi_v}(x)$ such that 
\eqref{E:nmin} holds. Let $v(\phi_{t^n}(x))=\alpha_n\in P_v$ and using the 
definition of $R_v(\alpha_n)$ for $\alpha_n\in P_v$ from \eqref{E:defR_v}, we 
get
\begin{equation}
\label{E:alpha_n}
|R_v((v(\phi_{t^n}(x)))|\le q^r.
\end{equation}
Inequality \eqref{E:alpha_n} and the result of lemma \eqref{L:L2} gives the 
estimate: 
\begin{equation}
\label{E:estimate}
|R_v(\alpha) |\le |R_v(v(\phi_{t^n}(x)))|\cdot q^{\frac{r^3-r^2}{2}}\le q^r\cdot 
q^{\frac{r^3-r^2}{2}}=q^{\frac{r^3-r^2+2r}{2}}
\end{equation}
for every $\alpha\in P_v(n)$ and for every $n\ge 0$.

Now, we estimate the magnitude of $n$, i.e. the number of steps that we need to 
make 
starting with $0<v(x)<N_v$ such that in the end $\phi_{t^n}(x)$ satisfies 
either \eqref{E:14} or \eqref{E:15}.

\begin{lemma}
\label{L:L13}
Assume $v\in S$ and $\hhat_v(x)>0$. Then there exists a set $P$ of cardinality 
bounded in terms of $r$ and $e(v|v_0)$ such that either 
$(v(x),\ac_{\pi_v}(x))\in P\times R_v(v(x))$ or 
$\hhat_v(x)>\frac{c_1}{e(v|v_0)^{\frac{r}{r_0}-1}d}$ with $c_1>0$ depending only 
on $\phi$.
\end{lemma}

\begin{proof}
If \eqref{E:15} does not hold for $x$ 
then we know that there exists $i_0\ge r_0$ such that 
$v(\phi_t(x))=q^{i_0}v(x)+v(a_{i_0})$.

Now, if $\phi_t(x)$ also does not satisfy \eqref{E:15} then for some $i_1$
$$v(\phi_{t^2}(x))=q^{i_1}v(x)+v(a_{i_1})\le q^iv(\phi_t(x))+v(a_i)$$
for all $i\in\{r_0,\dots,r\}$. So, in particular 
\begin{equation}
\label{E:43}
v(\phi_{t^2}(x))\le q^{i_0}v(\phi_t(x))+v(a_{i_0})
\end{equation}
and in general
\begin{equation}
\label{E:general}
v(\phi_{t^{k+1}}(x))\le q^{i_0}v(\phi_{t^k}(x))+v(a_{i_0})
\end{equation}
if $(v(\phi_{t^k}(x)),\ac_{\pi_v}(\phi_{t^k}(x)))\notin P_v\times 
R_v(v(\phi_{t^k}(x)))$.
Let us define the following sequence $(y_j)_{j\ge 0}$ by
$$y_0=v(x) \text{ and for all }j\ge 1 \text{: } 
y_j=q^{i_0}y_{j-1}+v(a_{i_0}).$$
If $\phi_{t^i}(x)$ does not satisfy \eqref{E:15} for $i\in\{0,\dots,n-1\}$ then 
by \eqref{E:general},
\begin{equation}
\label{E:22}
y_n\ge v(\phi_{t^n}(x)).
\end{equation}
The sequence $(y_j)_{j\ge 0}$ can be easily computed and we see that
\begin{equation}
\label{E:42}
y_j=q^{i_0j}\left(v(x)+\frac{v(a_{i_0})}{q^{i_0}-1}\right)-\frac{v(a_{i_0})}{q^{
i_0}-1}.
\end{equation}
But $v(x)<-\frac{v(a_{i_0})}{q^{i_0} -1}$, as a consequence of $v(x)<N_v$ and 
the proof of lemma \eqref{L:L1} (see equation \eqref{E:vi_0}). Thus,
\begin{equation}
\label{E:inegalitate1}
v(x)+\frac{v(a_{i_0})}{q^{i_0}-1}\le -\frac{1}{q^{i_0}-1}
\end{equation}
because $v(x),v(a_{i_0})\in \mathbb{Z}$.
Using inequality \eqref{E:inegalitate1} in the formula \eqref{E:42} we get
\begin{equation}
\label{E:inegalitate2}
y_j\le \frac{1}{q^{i_0}-1}(-q^{i_0j}-v(a_{i_0})).
\end{equation}
For $v_0$ the valuation of $K$ that sits under the valuation $v$ of $L$, we 
define
\begin{equation}
\label{E:cv_0}
c_{v_0}=\max\left\{ -v_0(a_i) \vert r_0\le i\le r\right\} .
\end{equation}
So, $c_{v_0}\ge 1$ because we know that at least one of the $a_i$ has a pole at 
$v$, 
thus at $v_0$ (we are working under the assumption that $v\in S$). Clearly, 
$c_{v_0}$ depends only on $\phi$ and on $K$; thus, for 
simplicity we denote $c_{v_0}$ by $c$ in the next calculations. Because of the 
definition of $c$, we have 
\begin{equation}
\label{E:inegalitate3}
-v(a_{i_0})\le e(v\vert v_0)c
\end{equation}
where $e(v\vert v_0)$ is as always the ramification index of $v$ over $v_0$.
Now, if we pick $m$ minimal such that
\begin{equation}
\label{E:23}
q^{r_0m}\ge ce(v\vert v_0)
\end{equation}
then we see that $m$ depends only on $\phi$ and $e(v\vert v_0)$. Using that 
$i_0\ge r_0$ we get that 
$$q^{i_0m}\ge ce(v\vert v_0).$$
So, using inequalities \eqref{E:inegalitate2}, \eqref{E:inegalitate3} and 
\eqref{E:23} we obtain $y_m\le 0$. Because of \eqref{E:22} we derive that
$$v(\phi_{t^m}(x))\le 0$$
which according to the dichotomy from lemma \eqref{L:L6'} yields that 
$\phi_{t^m}(x)$ satisfies either \eqref{E:14} or 
\eqref{E:15}. Thus, we need at most $m$ steps to get from $x$ to some 
$\phi_{t^n}(x)$ for which one of the two equations \eqref{E:14} or \eqref{E:15} 
is valid. This means that either 
\begin{equation}
\label{E:24}
\hat{h}_v(x)\ge \frac{e(v\vert v_0)}{q^{rm}q^{2r}d}\text{ (which holds if 
\eqref{E:14} is valid after $n\le m$ steps),}
\end{equation}
or
\begin{equation}
\label{E:25}
\phi_{t^n}(x) \, \text{satisfies \eqref{E:15} for} \, n\le m.
\end{equation}
This last equation implies that $(v(x),\ac_{\pi_v}(x))\in P_v(n)\times 
R_v(v(x))$ 
for some $n\le m$.

We analyze now the inequality from equation \eqref{E:24}. By the minimality of 
$m$ satisfying \eqref{E:23}, we have 
\begin{equation}
\label{E:26}
q^{rm}= (q^{r_0(m-1)})^{\frac{r}{r_0}} q^r< (ce(v\vert v_0) 
^{\frac{r}{r_0}}q^r.
\end{equation}
So, if \eqref{E:24} holds, we have the following inequality
\begin{equation}
\label{E:27}
\hat{h}_v(x)>\frac{e(v\vert v_0)}{c^{\frac{r}{r_0}}q^{3r}e(v\vert 
v_0)^{\frac{r}{r_0}}d}.
\end{equation}
We denote by $P=\bigcup_{i=0}^{m} P_v(i)$. We proved that for $i\ge 1$, $\vert 
P_v(i)\vert\le \vert P_v(0)\vert$ (and $P_v=P_v(0)$ has cardinality depending 
only on $r$; this was mainly the content of \eqref{L:L2}). To 
simplify the notations in the future we introduce new constants $c_i$, that will 
always depend only on $\phi$. For example, \eqref{E:27} says that
\begin{equation}
\label{E:28}
\hat{h}_v(x)>\frac{c_1}{e(v\vert v_0)^{\frac{r}{r_0}-1}d}  \text{ or }  
(v(x),\ac_{\pi_v}(x))\in P\times R_v(v(x))
\end{equation}
and $\vert R_v(v(x))\vert \le q^{\frac{r^3-r^2+2r}{2}}$ (see equation 
\eqref{E:estimate}), while $\vert P\vert\le \frac{r^2-r+2}{2}(m+1)$ with $m$ 
satisfying \eqref{E:26}.
\end{proof}

\section{The local Lehmer inequality}

We continue with the notation from Section $2$. The field $L$ is finitely 
generated and $v\in M_L$. First we will prove the following general lemma on 
valuations.

\begin{lemma}
\label{L:L-1}
Let $I$ be a finite set of integers. Let $N$ be an integer greater or equal than 
all the elements of $I$. For each $\alpha\in I$, let $R(\alpha)$ be a nonempty 
finite set 
of nonzero elements of the residue field at $v$. Let $W$ be an 
$\mathbb{F}_q$-vector subspace of $L$ with the property that for all $0\ne w\in 
W$, $(v(w),\ac_{\pi_v}(w))\in I\times R(v(w))$ whenever $v(w)\le N$.

Let $f$ be the smallest integer greater or equal than $\max_{\alpha\in 
I}\log_{q}\vert R(\alpha)\vert$. Then the codimension of $\left\{w\in W\mid 
v(w)>N\right\}$ is bounded by $\vert I\vert f$.
\end{lemma}
\begin{proof}
To prove \eqref{L:L-1} it suffices to show the following statement.
\begin{Claim}
\label{S:statement}
We cannot find a subspace $W'\subset W$ of dimension $1+\vert I\vert f$ such 
that for all $0\ne w\in W'$, $(v(w),\ac_{\pi_v}(w))\in I\times R(v(w))$.
\end{Claim}
To prove claim \eqref{S:statement} we will use induction on $z=\vert 
I\vert$.

If $z=1$, then $I=\{\alpha\}$. Assume that there are $(1+f)$ 
$\mathbb{F}_q$-linearly independent elements 
$$w_1,w_2,\dots,w_{f+1}$$
such that for all $0\ne w\in\Span(\{w_1,\dots,w_{f+1}\})$, $v(w)=\alpha$ and 
$\ac_{\pi_v}(w)\in R(\alpha)$. By our choice for $f$, we have more nonzero 
$\mathbb{F}_q$-linear combinations of 
$$\ac_{\pi_v}(w_1),\dots,\ac_{\pi_v}(w_{f+1})$$
than elements of $R(\alpha)$. Thus, there exists an $\mathbb{F}_q$-linear 
combination
$$\gamma=d_1\ac_{\pi_v}(w_1)+\dots+d_{f+1}\ac_{\pi_v}(w_{f+1})$$
where not all $d_1,\dots,d_{f+1}$ are $0$ and, either $0\ne\gamma\notin 
R(\alpha)$ or $\gamma=0$. So, if we let 
$$w=d_1w_1+\dots+d_{f+1}w_{f+1}$$
then, either $\ac_{\pi_v}(w)\notin R(\alpha)$ or $v(w)>\alpha$. Thus, claim 
\eqref{S:statement} holds for $z=1$.

Now we prove the inductive step. We assume that our claim \eqref{S:statement} 
holds for $z$, with $z\ge 1$, and we prove it for $(z+1)$. 

We assume $I=\{\alpha_1,\alpha_2,\dots,\alpha_{z+1}\}$ and we suppose there 
exist $\left(1+(z+1)f\right)$ $\mathbb{F}_q$-linearly independent elements of 
$W$, 
$w_1,\dots,w_{1+(z+1)f}$ such that for all nonzero 
$$w\in\Span(\{w_1,\dots,w_{1+(z+1)f}\})$$
we have that $(v(w),\ac_{\pi_v}(w))\in I\times R(v(w))$.

By the induction hypothesis we know that there are no $(1+zf)$ indices 
$$i_j\in\{1,\dots,1+(z+1)f\}$$ 
such that for each such index, $v(w_{i_j})\ge \alpha_2$. Thus, without loss of 
generality we may assume that there exists $1\le g\le 1+zf$ such that 
$v(w_i)=\alpha_1$ for all $i\in\{1,\dots,f+g\}$, while $v(\alpha_i)\ge \alpha_2$ 
if $i>f+g$. If $g=1+zf$, then there are no indices $i$ such that $v(w_i)\ge 
\alpha_2$.

There are two cases:

(i)  
$\dim_{\mathbb{F}_q}\Span(\{\ac_{\pi_v}(w_1),\dots,\ac_{\pi_v}(w_{f+g})\})>f$.

In this case, the definition of $f$ from \eqref{L:L-1} yields the existence of a 
nonzero 
$$w\in\Span(\{w_1,\dots,w_{f+g}\})$$
such that $v(w)=\alpha_1$ and $\ac_{\pi_v}(w)\notin R(\alpha_1)$. This proves 
the claim \eqref{S:statement} if (i) holds.

(ii)  
$\dim_{\mathbb{F}_q}\Span(\{\ac_{\pi_v}(w_1),\dots,\ac_{\pi_v}(w_{f+g})\})\le 
f$.

In this case, without loss of generality we may assume that 
$$\Span(\{ac_{\pi_v}(w_1),\dots,\ac_{\pi_v}(w_{f+g})\})=\Span(\{ac_{\pi_v}(w_1),
\dots,\ac_{\pi_v}(w_f)\}).$$ 
Then for every $k\in\{1,\dots,g\}$ we can find $l_k\in\Span(\{w_1,\dots,w_f\})$ 
such that $v(l_k)=\alpha_1$ and $\ac_{\pi_v}(l_k)=\ac_{\pi_v}(w_{f+k})$. Thus 
\begin{equation}
\label{E:-10}
v(w_{f+k}-l_k)>\alpha_1
\end{equation}
Because $0\ne w_{f+k}-l_k\in\Span(\{w_1,\dots,w_{1+(z+1)f}\})$ for every 
$k\in\{1,\dots,g\}$, it means that
\begin{equation}
\label{E:-20}
v(w_{f+k}-l_k)\in I.
\end{equation}
Equations \eqref{E:-10} and \eqref{E:-20} yield that $v(w_{f+k}-l_k)\ge 
\alpha_2$. This shows that for all $0\ne 
w\in\Span(\left\{w_{f+1}-l_1,\dots,w_{f+g}-l_g,w_{f+g+1},\dots,w_{1+(z+1)f} 
\right\})$, 
$$v(w)\in\{\alpha_2,\dots,\alpha_{z+1}\}\text{ and }\ac_{\pi_v}(w)\in R(v(w)).$$ 
Because all of the $(1+zf)$ elements 
$w_{f+1}-l_1,\dots,w_{f+g}-l_g,w_{f+g+1},\dots,w_{1+(z+1)f}$ are 
$\mathbb{F}_q$-linearly independent, we get a contradiction with the inductive 
hypothesis of claim \eqref{S:statement}. This concludes the proof of lemma 
\eqref{L:L-1}. 
\end{proof}

We are ready to prove theorem \eqref{T:T1}.
\begin{proof}[Proof of Theorem ~\ref{T:T1}]
First we observe that if $v\notin S$ then by lemma \eqref{L:L1'} we 
automatically get the lower bound $\hhat_v(x)\ge\frac{1}{d}$ because it must be 
that 
$v(x)<0$, otherwise we would have $\hhat_v(x)=0$. So, from now on we suppose 
that the valuation $v$ is from $S$.

We denote by $z=\vert P\vert$. Let $f$ be the smallest integer such that 
$$f\ge\max_{\alpha\in P}\log _q |R_v(\alpha)|.$$ 
So $f\le \frac{r^3-r^2+2r}{2}$, as shown by the proof of lemma \eqref{L:L13}. We 
also have the following inequality
\begin{equation}
\label{E:zf}
zf\le \frac{r^2-r+2}{2}(m+1)\cdot\frac{r^3-r^2+2r}{2}= 
\frac{r^5-2r^4+5r^3-4r^2+4r}{4}(m+1).
\end{equation}

Let $W=\Span(\left\{x,\phi_t(x),\dots,\phi_{t^{zf}}(x)\right\})$. Because 
$\hhat_v(x)>0$ we know that $x\notin\phi_{\tor}$ and so, 
$\dim_{\mathbb{F}_q}W=1+zf$. We also get from $\hhat_v(x)>0$ that for all $0\ne 
w\in W$, $\hhat_v(w)>0$. Then by lemma \eqref{L:L11}, we get that for all $0\ne 
w\in W$, $v(w)\le N_v-1$.

We apply lemma \eqref{L:L-1} to $W$ with $I=P$, $R=R_v$, $N=N_v-1$ and conclude 
that there exists $0\ne b\in \mathbb{F}_q[t]$, of degree at most $zf$ in $t$ 
such that 
\begin{equation}
\label{E:b}
(v(\phi_b(x)),\ac_{\pi_v}(\phi_b(x)))\notin P\times R_v(v(\phi_b(x))).
\end{equation}

We know that $\hat{h}_v(x)>0$ and so $\hat{h}_v(\phi_b(x))>0$. Equations 
\eqref{E:b} and \eqref{E:28} yield
$$\hat{h}_v(\phi_b(x))>\frac{c_1}{e(v\vert v_0)^{\frac{r}{r_0}-1}d}.$$
Thus
$$\hat{h}_v(x)>\frac{c_1}{q^{r\deg(b)}e(v\vert v_0)^{\frac{r}{r_0}-1}d}.$$
But, using inequality \eqref{E:zf}, we obtain 
$$q^{r\deg(b)}\le q^{rzf}\le 
q^{\frac{r^5-2r^4+5r^3-4r^2+4r}{4}r(m+1)}=q^{\frac{r^6-2r^5+5r^4-4r^3+4r^2}{4}
}(q^{rm})^{\frac{r^5-2r^4+5r^3-4r^2+4r}{4}}.$$
We use \eqref{E:26} and we get
$$q^{r\deg(b)}< q^{\frac{r^6-2r^5+5r^4-4r^3+4r^2}{4}}(ce(v\vert 
v_0))^{\frac{r}{r_0}\cdot\frac{r^5-2r^4+5r^3-4r^2+4r}{4}}q^{\frac{r^6-2r^5+5r^4-
4r^3+4r^2}{4}}.$$
Thus there exists a constant $C>0$ depending only on $c_1$, $c$, $q$ and $r$ 
such that 
\begin{equation}
\label{E:e}
\hat{h}_v(x)>\frac{C}{e(v\vert 
v_0)^{\frac{r}{r_0}(\frac{r^5-2r^4+5r^3-4r^2+4r}{4}+1)-1}d}.
\end{equation}
Because $c_1$ and $c$ depend only on $\phi$ we get the conclusion of 
\eqref{T:T1}.
\end{proof}

Using that $e(v\vert v_0)\le d$, we get the conclusion of theorem \eqref{T:T1'} 
$$\hat{h}(x)\ge\frac{C}{d^{k}}\text{ ,}$$
with $k\le\frac{r^6-2r^5+5r^4-4r^3+4r^2+4r}{4r_0}$.

\begin{remark}
\label{R:dependence}
From the above proof we see that the constant $C$ depends only on $q$, $r$ and 
the numbers $v(a_i)$ for $r_0\le i\le r-1$, in the hypothesis that $\phi_t$ is 
monic as a polynommial in $\tau$. As we said before, for the general case, when 
$\phi_t$ is not neccessarily monic, the constant $C$ from \eqref{T:T1} will be 
multiplied by the 
inverse of the degree of the extension of $K$ that we have to allow in order to 
construct a 
conjugated Drinfeld module $\phi^{(\gamma)}$ for which $\phi^{(\gamma)}_t$ is 
monic. The degree of this extension is at most $(q^r-1)$ because 
$\gamma^{q^r-1}a_r=1$.
\end{remark}
\begin{remark}
\label{R:R1}
It is interesting to note that \eqref{E:e} shows that the original statement of 
\eqref{C:Con2} holds, i.e. $k=1$, in the case that $e(v\vert v_0)=1$, which is 
the case when $x$ belongs to an unramified extension above $v_0$. Also, as 
observed in the beginning of the proof of \eqref{T:T1}, if $v$ and so, 
equivalently $v_0$ is not a pole for any of the $a_i$ then we automatically get 
exponent $k=1$ in theorem \eqref{T:T1}, as proved in lemma \eqref{L:L1'}. 

So, we see that in the course of proving \eqref{T:T1} we got an even stronger 
result that allows us to conclude that conjecture \eqref{C:Con2} and so, 
implicitly conjecture \eqref{C:Con1} holds in the maximal unramified extension  
above the finitely many irreducible divisors from $S_0$.
\end{remark}

\begin{remark}
\label{R:R2}
Also, it is interesting to note that the above proof shows that for every 
divisor $v$ associated to $L$ (as in Section $2$), there exists a number $n$ 
depending only on $r$ and $e(v\vert v_0)$ so that there exists 
$b\in\mathbb{F}_q[t]$ of degree at most $n$ in $t$ for which either 
$v(\phi_{b}(x))<M_v$ (in which case $\hhat_v(x)>0$), or $v(\phi_{b}(x))\ge N_v$ 
(in which case $\hhat_v(x)=0$). 
\end{remark}

\begin{example}
\label{E:E1}
  The result of theorem \eqref{T:T1} is optimal in the sense that we 
cannot hope to get the conjectured Lehmer inequality for the local height, 
i.e. $\frac{C}{d}$. We can only get, in the general case for the local height, 
an inequality with some exponent $k>1$, i.e. $\frac{C}{d^k}$.

For example, take $A=\mathbb{F}_q[t]$ and define 
$$\phi_t=\tau^r-t^{1-q}\tau .$$
Let $K=\mathbb{F}_q(t)$. Let $d=q^{m}-1$, for some $m\ge r$. Then let 
$x=t\alpha$ where 
$\alpha$ is a root of 
$$\alpha^{d}-\alpha-\frac{1}{t}=0.$$
Then $L=K(x)$ is totally ramified above $t$ of degree $d$. Let $v$ be the unique 
valuation of $L$ for which $v(t)=d$. We compute 
$$P_v=\left\{\frac{-d(q-1)}{q^r-q}\right\}$$
$$M_v=-\frac{d(q-1)}{q^r-q}$$
$$N_v=d$$
$$v(x)=d-1=q^m-2.$$
We compute easily $v(\phi_{t^i}(x))=d-q^i$ for every 
$i\in\{0,\dots,m\}$. 

Furthermore, $v(\phi_{t^m}(x))=d-q^m=-1\ne\frac{-d(q-1)}{q^r-q}$, 
because $\frac{-d(q-1)}{q^r-q}\notin\mathbb{Z}$ ($q\not\vert d(q-1)$). Thus 
$v(\phi_{t^m}(x))$ is negative and not in $P_v$ and so, \eqref{L:L3'} yields
$$v(\phi_{t^{m+1}}(x))<M_v.$$
Actually, because $m\ge r$, an easy computation shows that  
$$v(\frac{\phi_{t^m}(x)^q}{t^{q-1}})=-q-d(q-1)=-q^{m+1}+q^m-1<-q^r= 
v((\phi_{t^m}(x))^{q^r}).$$
This shows that $v(\phi_{t^{m+1}}(x))=-q^{m+1}+q^m-1<M_v<0$ and so, by 
\eqref{L:L2'} 
$$\hhat_v(x)=\frac{\hhat_v(\phi_{t^{m+1}}(x))}{q^{r(m+1)}}=\frac{q^{m+1}-q^m+1} 
{q^{r(m+1)}d}<\frac{q^{m+1}} 
{q^{m+r}q^{(r-1)m}d}<\frac{q^{1-r}}{d^r},$$
because $d=q^m-1<q^m$.

This computation shows that for Drinfeld modules of type  
$$\phi_t=\tau^r-t^{1-q}\tau$$
the exponent $k$ from \eqref{T:T1} should be at least $r$. The exact same 
computation will give us that in the case of a Drinfeld module of the form
$$\phi_t=\tau^r-t^{1-q^{r_0}}\tau^{r_0}$$
for some $1\le r_0<r$ and $x$ of valuation $\left(q^{r_0m}-2\right)$ at a place 
$v$ that 
is totally ramified above the place of $t$ with ramification index $q^{r_0m}-1$, 
the exponent 
$k$ in theorem \eqref{T:T1} should be at least 
$\frac{r}{r_0}$. In theorem \eqref{T:T2} we will prove that for non-wildly 
ramified extensions above places from $S_0$, we indeed get exponent 
$k=\frac{r}{r_0}$. But before doing this, we observe that the present example is 
just a counter-example to statement \eqref{C:Con2}, not to conjecture 
\eqref{C:Con1}. In other words, the global Lehmer inequality holds for our 
example even if the local one fails. 

Indeed, because $x$ was chosen to have positive valuation at the only place from 
$S$, it means that there exists another place, call it $v'$ which is not in $S$, 
for which $v'(x)<0$. But then by lemma \eqref{L:L1'}, we get that 
$\hhat_{v'}(x)\ge\frac{1}{d}$, which means that also $\hhat(x)\ge\frac{1}{d}$. 
Thus we obtain a lower bound for the global height as conjectured in 
\eqref{C:Con1}.
\end{example}

Now, in order to get to the result of \eqref{T:T2} we prove a lemma.
\begin{lemma}
\label{L:lemma3}
With the notation from the proof of theorem \eqref{T:T1}, let 
$L=\lcm_{i\in\{1,\dots,r-r_0\}}\lbrace q^i-1\rbrace$. If $p$ does not divide 
$e(v\vert v_0)$ then $e(v\vert v_0)$ divides $L\alpha$ for every $\alpha\in P$.
\end{lemma}
\begin{proof}  
Indeed, from its definition \eqref{E:defP_v}, $P_v$ contained $\{0\}$ and 
numbers of the form 
$$\frac{v(a_i)-v(a_j)}{q^j-q^i}=\frac{v(a_i)-v(a_j)}{q^i(q^{j-i}-1)},$$
for $j>i$.
Clearly, every number of this form, times $L$ is divisible by $e(v\vert v_0)$, 
because we supposed that $p\not\vert e(v\vert v_0)$. The set $P_v(1)$ contains 
numbers of the form
\begin{equation}
\label{E:process}
\frac{\alpha -v(a_i)}{q^i}
\end{equation}
where $\alpha\in P_v=P_v(0)$ and $a_i\ne 0$. Using again that $p$ does not 
divide $e(v|v_0)$ we get that $e(v|v_0)\mid L\alpha_1$ for all $\alpha_1\in 
P_v(1)$. Repeating the process from \eqref{E:process} we obtain all the elements 
of $P_v(n)$ for every $n\ge 1$ and by induction on $n$, we conclude that 
$e(v|v_0)\mid L\alpha_n$ for all $\alpha_n\in P_v(n)$. Because 
$P=\bigcup_{n=0}^{m}P_v(n)$ we get the result of this lemma.
\end{proof}

\begin{theorem}
\label{T:T2}
Let $\phi:A\rightarrow K\lbrace\tau\rbrace$ be a Drinfeld module of finite 
characteristic. Let $t\in A$ such that $\phi_t=\sum_{i=1}^r a_i\tau^i$ is 
inseparable. Let $r_0$ the index of the first nonzero coefficient of 
$\phi_t$. Let $x\in K^{\alg}$ and let $v\in M_{K(x)}$ such that $h_v(x)>0$. Let 
$v_0$ be the valuation on $K$ that sits below $v$. 

If $p$ does not divide $e(v|v_0)$, there exists a constant $C>0$ depending only 
on $\phi$ such that 
$\hat{h}_v(x)\ge\frac{C}{e(v|v_0)^{\frac{r}{r_0}-1}[K(x):K]}$.
\end{theorem}
\begin{proof} 

Just as we observed in Section $2$ and in remark \eqref{R:dependence}, it 
suffices to prove \eqref{T:T2} under the hypothesis that $\phi_t$ is monic in 
$\tau$.
 
Let now $d=[K(x):K]$. We observe again that from \eqref{L:L1'} it follows that 
if $v\notin S$ then 
$\hhat_v(x)\ge\frac{1}{d}\ge\frac{1}{e(v|v_0)^{\frac{r}{r_0}-1}d}$. 
So, from now on we consider the case $v\in S$.

Then, using the result of \eqref{L:lemma3} in \eqref{E:vi_0} we see that 
\begin{equation}
\label{E:inegalitate1'}
v(x)+\frac{v(a_{i_0})}{q^{i_0}-1}\le -\frac{\frac{e(v\vert v_0)}{L}}{q^{i_0}-1}
\end{equation}
if $v(x)\in P$. Then also \eqref{E:inegalitate2} changes into
\begin{equation}
\label{E:inegalitate2'}
x_m\le \frac{1}{q^{i_0}-1}(-q^{i_0m}\frac{e(v\vert w)}{L} -v(a_{i_0})).
\end{equation}
So, then we choose $m'$ minimal such that
\begin{equation}
\label{E:inegalitate3'}
q^{r_0m'}\ge cL
\end{equation}
where $c=c_{v_0}$ is the same as in \eqref{E:cv_0}. Thus $m'$ depends only on 
$\phi$. We redo the computations from \eqref{E:24} to 
\eqref{E:28}, this time with $m'$ in place of $m$ and because of 
\eqref{E:inegalitate2'} and \eqref{E:inegalitate3'}, we get that 
\begin{equation}
\label{E:29}
\hat{h}_v(x)>\frac{c_1}{e(v\vert v_0)^{\frac{r}{r_0}-1}d}  \text{ or }  
(v(x),\ac_{\pi_v}(x))\in P'\times R_v(v(x))
\end{equation}
where $P'=\bigcup_{i=0}^{m'}P_v(i)$. At this moment we can redo the argument 
from the proof of \eqref{T:T1} using $P'$ instead of $P$, only that now 
$z'=\vert P'\vert$ is independent of $x$ or $d$. We conclude once again that 
there exists $b$, a polynomial in $t$ of degree at most $z'f$ such that
$$\hat{h}_v(\phi_b(x))>\frac{c_1}{e(v|v_0)^{\frac{r}{r_0}-1}d}.$$
But because both $f$ and $z'$ depend only on $\phi$, we conclude that indeed, 
$$\hat{h}(x)\ge\frac{C}{e(v|v_0)^{\frac{r}{r_0}-1}d}$$
with $C>0$ depending only on $\phi$.
\end{proof}

\begin{example}
\label{E:E2}
We discuss now the conjecture \eqref{C:Con2} for Drinfeld modules of generic 
characteristic. So, consider the Carlitz defined on $\mathbb{F}_p[t]$ by 
$\phi_t=t\tau^0+\tau$, where $\tau(x)=x^p$ for all $x$. Take 
$K=\mathbb{F}_p(t)$. Let $L$ be a finite extension of $K$ which is totally 
ramified above $\infty$ 
and so, let the ramification index equals $d=[L:K]$. Also, let $v$ be the unique 
valuation of $L$ sitting above $\infty$.

Let $x\in L$ of valuation $nd$ at $v$ for some $n\ge 1$. An easy computation 
shows that for all $m\in\{1,\dots,n\}$, $v(\phi_{t^m}(x))=dn-dm$. So, in 
particular $v(\phi_{t^n}(x))=0$ and so, 
$$v(\phi_{t^{n+1}}(x))=-d<M_v=\frac{-d}{p-1}\text{ .}$$
This shows, after using lemma \eqref{L:L2'}, that 
$\hhat_v(\phi_{t^{n+1}}(x))=\frac{d}{d}=1$. This in turn implies that 
$$\hhat_v(x)=\frac{1}{p^{n+1}}.$$
But we can take $n$ arbitrarily large, which shows that there is no way to 
obtain a similar result like theorem \eqref{T:T1} for generic characteristic 
Drinfeld modules.
\end{example}

The next theorem shows that the example \eqref{E:E2} is in some sense the only 
way theorem \eqref{T:T1} fails for Drinfeld modules of generic characteristic.

\begin{theorem}
\label{T:T3}
Let $\phi$ be a Drinfeld module of generic characteristic and so, with the usual 
notation, let $\phi_t=t\tau^0+\sum_{i=1}^{r}a_i\tau^i$, for a non-constant 
$t\in A$. Let $x\in K^{\alg}$ and let $v$ be an irreducible divisor from  
$M_{K(x)}$ 
that does not sit over the place $\infty$ from $\Frac(A)$. Let $v_0\in M_K$ sit 
below $v$. There exist two 
positive constants $C$ and $k$ depending only on $\phi$, such that if 
$\hhat_v(x)>0$ then $\hhat_v(x)\ge\frac{C}{e(v|v_0)^{k-1}[K(x):K]}$.
\end{theorem}
\begin{proof} 

Again, as we mentioned in Section $2$ and in remark \eqref{R:dependence}, it 
suffices to prove this theorem under the hypothesis that $\phi_t$ is monic in 
$\tau$. Also, if $v\notin S$ theorem \eqref{T:T3} holds as shown by lemma 
\eqref{L:L1'}.
 
The analysis of local heights done in Section $2$ applies to both finite and 
generic characteristic until lemma \eqref{L:L11}. So, we still get the 
conclusion of lemma \eqref{L:L6'}. Thus, if $v(x)\le 0$ then either 
$\hhat_v(x)\ge\frac{e(v|v_0)}{q^{2r}[K(x):K]}$ or $(v(x),\ac_{\pi_v}(x))\in 
P_v\times R_v(v(x))$, with $\vert 
P_v\vert$ and $\vert R_v(v(x))\vert $ depending only on $q$ and $r$.

We know from our hypothesis that $v(t)\ge 0$ and so,
\begin{equation}
\label{E:41}
v(tx)\ge v(x).
\end{equation}
Now, if $v(x)\ge N_v$, then $v(a_ix^{q^i})\ge v(x)$, for all $i\ge 1$ (by the 
definition of $N_v$) and using also equation \eqref{E:41}, we get
$$v(\phi_t(x))\ge v(x)\ge N_v.$$
Iterating this computation we get that $v(\phi_{t^n}(x))\ge N_v$, for all $n\ge 
1$ and so $\hhat_v(x)=0$, contradicting the hypothesis of our theorem. This 
argument is the equivalent of lemma \eqref{L:L11} for Drinfeld modules of 
generic characteristic under the hypothesis $v(t)\ge 0$.

Thus it must be that $v(x)<N_v$. Then, lemma \eqref{L:L1} holds identically. 
This 
yields that either 
$(v(x),\ac_{\pi_v}(x))\in P_v\times R_v(v(x))$ or $v(\phi_t(x))<v(x)$.

From this point on, the proof continues just as for theorem \eqref{T:T1}. We 
form 
just as before the sets $P_v(n)$ and their union will be again denoted by $P$. 
We conclude once again as in \eqref{E:27} that \emph{either} 
$$\hhat_v(x)\ge\frac{1}{q^{3r}c^{\frac{r}{r_0}}e(v\vert 
v_0)^{\frac{r}{r_0}-1}[K(x):K]}$$ 
with the same $c>0$ depending only on $q$, 
$r$ and $\phi$ as in the proof of \eqref{T:T1}, \emph{or} 
$$(v(x),\ac_{\pi_v}(x))\in 
P\times R_v(v(x))$$
where $\vert P\vert$ is of the order of $\log e(v\vert v_0)$. We observe that 
when we use equations \eqref{E:general}, \eqref{E:42}, \eqref{E:inegalitate1}, 
\eqref{E:inegalitate2} the index $i_0$ is still at least $1$. This is the case 
because if 
$v(x)<N_v$ and $(v(x),\ac_{\pi_v}(x))\notin P_v\times R_v(v(x))$ then there 
exists 
$i_0\ge 1$ such that $v(\phi_t(x))=v(a_{i_0})+q^{i_0}v(x)$. Also, $v(x)<N_v$ 
means that there exists at least one index $i\in\{1,\dots,r\}$ such 
that $v(tx)\ge v(x)>v(a_ix^{q^i})$.

Thus, the first index other than $0$ of a non-zero coefficient $a_i$ will play 
the role of $r_0$ as in the proof of 
\eqref{T:T1}. Finally, lemma \eqref{L:L-1} finishes the proof of theorem 
\eqref{T:T3}.
\end{proof}

So, we get in the same way as in the proof of \eqref{T:T1}, the conclusion for 
theorem \eqref{T:T3}. The difference made by $v$ not sitting above $\infty$ is 
that for 
$v(x)\ge 0$, $v(\phi_t(x))$ can decrease only if $v(x)<N_v$, i.e. only if 
there exists $i\ge 1$ such that $v(a_ix^{q^i})<v(x)$. If $v$ sits over $\infty$, 
then 
$v(tx)<v(x)$ and so, $v(\phi_t(x))$ might decrease just because of the $t\tau^0$ 
term from $\phi_t$. Thus, in that case, as example \eqref{E:E2} showed, we can 
start with $x$ having arbitrarily large valuation and we are able to 
decrease it by applying $\phi_t$ to it repeatedly, making the valuation of 
$\phi_{t^n}(x)$ be less than $M_v$, which would mean that $\hhat_v(x)>0$. But in 
doing this we will 
need a number $n$ of steps (of applying $\phi_t$) that we will not be able to 
control; so $\hhat_v(x)$ will be arbitrarily small.

It is easy to see that remarks \eqref{R:R1} and \eqref{R:R2} are valid also for 
theorem \eqref{T:T3} in the hypothesis that $v$ does not sit over the place 
$\infty$ of $\Frac(A)$. Also, just as we were able to derive theorem 
\eqref{T:T2} from 
the proof of \eqref{T:T1}, we can do the same thing in theorem 
\eqref{T:T3} and find a specific value of the constant $k$ that will work in the 
case that $v$ is not wildly ramified above $v_0\in M_K$. The result is the 
following theorem whose proof goes along the same 
lines as the proof of \eqref{T:T2}.  
\begin{theorem}
\label{T:T4}
Let $\phi$ be a Drinfeld module of generic characteristic and let 
$\phi_t=t\tau^0+\sum_{i=r_0}^{r}a_i\tau^i$, with $a_{r_0}\ne 0$ (of course, 
$r_0\ge 1$). There exists a constant $C>0$, depending only on $\phi$ such that 
for every $x\in K^{\alg}$ and every $v\in M_{K(x)}$ that is neither wildly 
ramified above $K$ nor sitting above the place $\infty$ of $\Frac(A)$, if 
$\hhat_v(x)>0$ then 
$\hhat_v(x)\ge\frac{C}{e(v|v_0)^{\frac{r}{r_0}-1}[K(x):K]}$.
\end{theorem}

We can also construct an example simliar to \eqref{E:E1} which shows that 
constant $k=\frac{r}{r_0}$ is optimal in the above theorem. Indeed, if we take a 
Drinfeld module $\phi$ defined on $\mathbb{F}_q[t]$ by 
$$\phi_t=t\tau^0+\frac{1}{t^{q^{r_0}-1}}\tau^{r_0}+\tau^{r}$$
and $x$ as in example \eqref{E:E1} then a similar computation will show that we 
cannot hope for an exponent $k$ smaller than $\frac{r}{r_0}$.

The constants $C$ from theorem \eqref{T:T2}, \eqref{T:T3} and \eqref{T:T4} and 
the constant $k$ from \eqref{T:T3} have the same corresponding dependency on 
$q$, $r$ and $\phi$ as explained in the proof of theorem \eqref{T:T1}.

\address{Dragos Ghioca, Department of Mathematics, University of California, 
Berkeley, CA 94720}

\email{dghioca@math.berkeley.edu}


\begin{thebibliography}{9}
\bibitem{Den}
Laurent Denis, \emph{The Lehmer problem in finite characteristic}, (French) 
Compositio Math. \textbf{98} (1995), no. 2, 
167-175
\bibitem{DG}
Dragos Ghioca, \emph{Lehmer inequality and the Mordell-Weil theorem for 
Drinfeld modules}, preprint 2004
\bibitem{Goss}
David Goss, \emph{Basic structures of function field arithmetic}, Ergebnisse 
der Mathematik und ihrer Grenzgebiete (3) [Results in Mathematics and Related
  Areas (3)], \textbf{35}. Springer-Verlag, Berlin, 1996
\bibitem{Poo} 
Bjorn Poonen, \emph{Local height functions and the Mordell-Weil theorem 
for Drinfeld modules}, Compositio Mathematica \textbf{97} (1995), 349-368
\bibitem{Wan}
Julie Tzu-Yueh Wang, \emph{The Mordell-Weil theorems for Drinfeld modules over 
finitely generated function fields}, Manuscripta math. \textbf{106}, 305-314 
(2001)
\end{thebibliography}
\end{document}